\documentclass[12pt]{article}
          \usepackage{amsfonts}
          \usepackage{amssymb}

\newtheorem{thm}{Theorem}

\newtheorem{lemma}[thm]{Lemma}

\newtheorem{prop}[thm]{Proposition}

\def\pf{\medbreak\noindent{\bf Proof:}\enspace}

\def\be{\begin{eqnarray}}
\def\ee{\end{eqnarray}}
\def\bee{\begin{eqnarray*}}
\def\eee{\end{eqnarray*}}

\def\ts{\textstyle}
\def\bra{\langle}
\def\ket{\rangle}

\def\rt2{\ts \frac{1}{\sqrt{2}} }

\def\la{\langle}
\def\ra{\rangle}

            \title{On the existence of a common quadratic
Lyapunov function for a rank one difference}

         \author{Christopher King  and Michael
Nathanson
\\
\\ Department of Mathematics
\\ Northeastern University
\\ Boston MA 02115
\\
{\normalsize king@neu.edu,  nathanson.m@neu.edu}
}

\begin{document}

\maketitle

\begin{abstract}
Suppose that $A$ and $B$ are real stable matrices, and that their difference
$A-B$ is rank one. Then $A$ and $B$ have a common quadratic Lyapunov function if and
only if the product $AB$ has no real negative eigenvalue. This result is due
to Shorten and Narendra, who showed that it follows as a consequence
of the Kalman-Yacubovich-Popov solution of the Lur'e problem. Here we
present a new and independent proof based on results from convex analysis
and the theory of moments.
\end{abstract}

\bigskip
\par\noindent{\bf Key words:} quadratic Lyapunov function; Hankel matrix;
discrete moment problem.

\bigskip
\par\noindent{\bf AMS Subject Classification:} 15A48, 34D20, 34H05

\pagebreak



\section{Introduction and statement of results}
This paper presents a new proof of the Shorten-Narendra Theorem \cite{SN}, which gives a
simple spectral condition for the existence of a common quadratic Lyapunov function
(CQLF) for two stable matrices whose difference is rank one. Recall that a
matrix $A$ is {\it stable} if the spectrum of
$A$ lies wholly in the open left half of the complex plane. An equivalent condition
is the existence of a 
positive definite matrix $P$ such that $P A + A^{\rm T} P$ is negative definite,
in which case the function ${\bf x}^{\rm T} P {\bf x}$ is a quadratic 
Lyapunov function for the system $\dot{\bf x} = A {\bf x}$.
Consideration of the switching system $\dot{\bf x} = A(t) {\bf x}$ with
$A(t) \in \{A,B\}$ leads to the notion of a common quadratic Lyapunov function (CQLF),
which is determined by a positive definite matrix $P$ satisfying
\be\label{def:CQLF}
P A + A^{\rm T} P < 0, \quad
P B + B^{\rm T} P < 0
\ee

The following theorem of Shorten and Narendra provides a simple
test for the existence of a CQLF in the case where $A - B$ is rank one.
The theorem is stated in \cite{SN} for matrices in companion form,
however this is unnecessary \cite{Meyer} and we state the 
result in its full generality here.

\medskip

\begin{thm}\label{thmSN}{\rm [Shorten and Narendra]}
Let $A$ and $B$ be stable matrices and suppose that
$A - B$ is rank 1. Then the necessary and
sufficient condition that $(A,B)$ have a CQLF is that the matrix
$A B$ does not have a real negative eigenvalue.
\end{thm}

\medskip
The proof of Theorem \ref{thmSN} presented in \cite{SN} 
first relates the spectral condition on $A B$ to the following
positivity condition for the resolvent of $A$ along the imaginary axis:
\be\label{cc}
{\rm Re} z = 0 \,\, \Rightarrow
1 + {\rm Re} \,\, {\bf v}^{\rm T} (z - A)^{-1} {\bf u}  > 0
\ee
where $A - B = {\bf u} {\bf v}^{\rm T}$.
The authors then make use of
earlier work of Narendra and Goldwyn \cite{NG} and Willems \cite{Will} which showed
that this resolvent condition (known as the circle criterion) is equivalent to the
existence of a CQLF. These earlier papers proved the equivalence by transforming
the question of the existence of a CQLF for the pair $(A, B)$ into the existence of 
a solution for the Lur'e problem. They then used the fundamental
results of Kalman \cite{Kal} who 
used techniques from analytic function theory to show that
the circle criterion gives a necessary and sufficient condition for the existence of
a solution of the Lur'e problem.

\medskip
The result of Theorem \ref{thmSN} is strikingly simple, 
and it gives an easy way to check for the existence of a CQLF. It also encourages the
belief that there should be a direct and independent proof which does not use the equivalence
between the CQLF problem and the Lur'e problem.
In this paper we provide such a proof, using methods
of convex analysis and the theory of moments. The proof has a
geometrical flavor which is described in the next paragraph.

\medskip
There is a dual formulation of the CQLF condition in terms
of intersecting cones in the space of symmetric matrices. Given a real matrix $A$,
define
\be\label{def:cone}
{\cal C}(A) = \{ A X + X A^{\rm T} \, | \, X \geq 0 \}
\ee
That is, ${\cal C}(A)$ is the cone of symmetric matrices of the form
$A X + X A^{\rm T}$ where
$X$ runs over all positive semidefinite matrices. Considering real $n \times n$ 
matrices as $n^2$-component vectors with the Hilbert-Schmidt inner product,
the existence of a quadratic Lyapunov function for $A$ is equivalent to the existence of a
positive definite matrix $P$ such that $\la P, M \ra  = {\rm Tr} P M < 0$ for all $M \neq 0$ in
${\cal C}(A)$. It is convenient to view this in terms of the hyperplane
which is the orthogonal complement of $P$, in which case the condition is that the cone
${\cal C}(A)$ lies on one side of the hyperplane.
Correspondingly, the
existence of a CQLF for $A$ and $B$ is equivalent to finding
such a hyperplane with both cones  ${\cal C}(A)$ and ${\cal C}(B)$
on the same side, or alternatively with the cones ${\cal C}(A)$ and ${\cal C}(-B)$
on opposite sides. Therefore the existence of a CQLF for $A$ and $B$ is
equivalent to the non-intersection (except at the origin) of the cones
${\cal C}(A)$ and ${\cal C}(-B)$. This observation leads to the
following proposition.

\medskip

\begin{prop}\label{prop1}
Let $A$ and $B$ be stable matrices.  Then the pair $(A, B)$ does {\rm NOT} have
a CQLF if and only if there are nonzero positive semidefinite matrices 
$X$ and $Y$ such that
\be
A X + X A^{\rm T} + B Y + Y B^{\rm T} = 0
\ee
\end{prop}

\medskip
The main result of this paper is contained in the following theorem.
It describes a special property of the intersection of the cones
${\cal C}(A)$ and ${\cal C}(-B)$ in the case of interest here, namely
when $A - B$ is rank 1. The extreme points of the cone ${\cal C}(A)$
have the form $A {\bf v} {\bf v}^{\rm T} + {\bf v} {\bf v}^{\rm T} A^{\rm T}$,
where ${\bf v}$ is a vector. The next theorem shows that
whenever the cones ${\cal C}(A)$ and ${\cal C}(-B)$ have a 
nonzero intersection, then the extreme points of the cones must also
have a nonzero intersection.

\medskip
\begin{thm}\label{thm2}
Let $A$ and $B$ be stable matrices, with
$A - B$ rank one.
Suppose that there are nonzero positive semidefinite matrices 
$X$ and $Y$ such that
\be\label{thm2:1}
A X + X A^{\rm T} + B Y + Y B^{\rm T} = 0
\ee
Then there are nonzero vectors ${\bf v}$ and ${\bf w}$ such that
\be\label{thm2:2}
A {\bf v} {\bf v}^{\rm T} + {\bf v} {\bf v}^{\rm T} A^{\rm T}
+ B {\bf w} {\bf w}^{\rm T} + {\bf w} {\bf w}^{\rm T} B^{\rm T}
= 0
\ee
\end{thm}

\medskip
Combining Proposition \ref{prop1} and Theorem \ref{thm2}
shows that the pair $(A, B)$ does {\it not} have a CQLF if and only if 
there are nonzero vectors ${\bf v}$ and ${\bf w}$ such that
(\ref{thm2:2}) holds.
The proof of Theorem \ref{thmSN} is completed by showing that for stable matrices $A$ and $B$,
the existence of vectors ${\bf v}$ and ${\bf w}$ satisfying (\ref{thm2:2})
is equivalent to the condition that $A B$ has a real negative eigenvalue.
This equivalence was first shown in a more general setting by Mason and Shorten \cite{MS}. 
The idea is simple: there are only two
possible ways for the equation (\ref{thm2:2}) to hold -- either 
${\bf v} = \alpha B {\bf w}$ for some $\alpha$, or else ${\bf v} = \alpha {\bf w}$
for some $\alpha$. The first possibility leads to
$(\alpha A B + \alpha^{-1}) {\bf w} = 0$, which is precisely the condition
that $A B$ has eigenvalue $- \alpha^{-2}$.  Running the argument in reverse shows that
the conditions are equivalent. 

The second possibility would imply that $(\alpha^2 A + B) {\bf w} = 0$,
or equivalently that some convex combination $(1-x) A + x B$ is singular.
However writing $A - B = R$ we have
\be\label{case2}
det[(1-x) A + x B] & = & det[A] \,\, det[I - x A^{-1} R] \nonumber \\
& = & det[A] \,\, \Big(1 - x Tr(A^{-1} R) \Big)
\ee
The left side of (\ref{case2}) is nonzero and has the same sign at $x=0$ and
$x=1$, hence the right side cannot vanish for any value of $x$ between $0$ and
$1$. This rules out the second possibility.

\medskip
Thus we see that Theorem \ref{thmSN} follows  directly from Theorem \ref{thm2},
and the rest of the paper is devoted to its proof.
In section 2 we show that it is sufficient to
assume a special form for the matrices known as companion form.  The main work of the paper
appears in section 3 where we prove Theorem \ref{thm2}.
The proof uses some relations between Hankel matrices and the
solution of the discrete moment problem, and some
linear algebra arguments. These are derived
in Appendices A and B.

\section{Reduction to companion form}
Let us write
\be
A - B = {\bf x} {\bf y}^{\rm T}
\ee
where ${\bf x}$ and $\bf y$ are vectors in ${\bf R}^n$. Let $V$ be the span of
${\bf x}, A {\bf x}, A^2 {\bf x}, \dots$. Suppose first that $V$ is a proper
subspace of ${\bf R}^n$.
Then ${\bf R}^n = V \oplus V^{\perp}$
and with respect to this decomposition $A$ and $B$ are block matrices of the form
\be\label{block1}
A = \pmatrix{A_1 & A_2 \cr 0 & A_3}, \quad\quad
B = \pmatrix{B_1 & B_2 \cr 0 & A_3}
\ee
Since the spectrum of $A$ is the union of the spectra of $A_1$ and $A_3$,
it follows that $A_1$, $A_3$ and $B_1$ are also stable. Now
suppose that (\ref{thm2:1}) holds, and write $X$ and $Y$ in block form
\be\label{block2}
X = \pmatrix{X_1 & X_2 \cr X_2^{\rm T} & X_3}, \quad\quad
Y = \pmatrix{Y_1 & Y_2 \cr Y_2^{\rm T} & Y_3}
\ee
Then it follows from (\ref{thm2:1}) that
\be
A_3 (X_3 + Y_3) + (X_3 + Y_3) A_3^{\rm T} = 0
\ee
which in turn implies that $X_3 + Y_3 = 0$ since $A_3$ is stable.
Then the positivity of $X$ and $Y$ imply that $X_2 = X_3 = Y_2 = Y_3 = 0$.
Therefore (\ref{thm2:1}) reduces to
\be\label{reduce}
A_1 X_1 + X_1 A_1^{\rm T} + B_1 Y_1 + Y_1 B_1^{\rm T} = 0
\ee
This means that it is sufficient to prove Theorem \ref{thm2} for the pair
$(A_1, B_1)$. Since $A_1 - B_1 = {\bf x} {\tilde {\bf y}}^{\rm T}$
where ${\tilde {\bf y}}$ is the projection of ${\bf y}$ onto $V$, 
the equation (\ref{reduce}) is a special case of (\ref{thm2:1}),
namely the case where the vectors ${\bf x}, A {\bf x}, A^2 {\bf x}, \dots$ span
the whole space. Hence without loss of generality we will assume
that the vectors ${\bf x}, A {\bf x}, A^2 {\bf x}, \dots, A^{n-1} {\bf x}$
are linearly independent. In this case the pair $(A, {\bf x})$ is called
{\it completely controllable}. We next show how this allows a change of basis into
a special form known as {\it companion form} (see \cite{Lef} for details).

We first introduce the matrix
\be
S = \pmatrix{0 & 1 & 0 & \dots & 0 \cr
0 & 0 & 1 & \dots & 0 \cr
\vdots && \ddots && \vdots \cr
0 && \dots & 0 & 1 \cr
0 && \dots & 0 & 0}
\ee
and the vector
\be
{\bf g} = \pmatrix{0 \cr \vdots \cr 0 \cr 1}
\ee
Then a matrix $A$ is said to be in companion form if it can be written
\be\label{def:comp}
A = S + {\bf g} {\bf h}^{\rm T}
\ee
for some vector $\bf h$. 
Suppose that  the pair $(A, {\bf x})$ is completely controllable,
and that the characteristic polynomial for $A$ is
\be\label{charA}
A^n + a_n A^{n-1} + \cdots + a_1 I = 0
\ee 
Then we choose the following vectors as a basis:
\be
{\bf e}_n & = & {\bf x} \nonumber  \\
{\bf e}_{n-1} & = & (A + a_n I ){\bf x} \nonumber  \\
{\bf e}_{n-2} & = & (A^2 + a_n A + a_{n-1} I){\bf x} \nonumber  \\
\vdots && \nonumber \\
{\bf e}_1 & = & (A^{n-1} + a_{n} A^{n-2} + \cdots + a_2 I){\bf x} \nonumber
\ee
The condition that $(A, {\bf x})$ is completely controllable guarantees
that these vectors form a  basis. Furthermore when the matrix $A$ is written
in this basis, it is easily seen to have the form (\ref{def:comp}).
That is, there is a non-singular real matrix $R$ such that
$R {\bf x} = {\bf g}$ and $R A R^{-1}$ has the form (\ref{def:comp}).
It follows that $R B R^{-1} = R A R^{-1} - {\bf g} {\bf y}^{\rm T} R^{-1}$ 
is also in companion form, and
furthermore that the pair $(A, B)$ has a CQLF  if and only if the pair
$(R A R^{-1}, R B R^{-1})$ has a CQLF. Similarly for the
condition that  $A B$ has a negative real eigenvalue.
Therefore it is sufficient to prove Theorem \ref{thmSN} for the
case that both $A$ and $B$ are in companion form.

\section{Proof of Theorem \ref{thm2}}
\subsection{Solving $A X + X A^{\rm T} + B Y + Y B^{\rm T} = 0$}
We assume henceforth that $A$ and $B$ are both in companion form, that is
\be\label{comp-form}
A = S + {\bf g} {\bf h}^{\rm T}, \quad
B = S + {\bf g} {\bf k}^{\rm T}
\ee
and that equation (\ref{thm2:1}) holds. Define
\be\label{def:Z}
Z = X + Y
\ee
and 
\be\label{def:w}
{\bf w} = {\bra {\bf k-h}, Y ({\bf k-h}) \ket}^{-1/2} \,\, Y ({\bf k-h})
\ee
Note that $\bra {\bf k-h}, Y ({\bf k-h}) \ket$ cannot be zero, as this would imply
that $A Z + Z A^{\rm T} = 0$ which is impossible since $Z \neq 0$.
Then we can rewrite (\ref{thm2:1}) as
\be\label{start2}
SZ + Z S^{\rm T} + {\bf g} {\bf h}^{\rm T} Z +
Z {\bf h} {\bf g}^{\rm T} + 
\bra {\bf k-h}, {\bf w} \ket \, \bigg(
{\bf g} {\bf w}^{\rm T}
+ {\bf w} {\bf g}^{\rm T} \bigg) = 0
\ee

It will be convenient to separate (\ref{start2}) into a pair of equations.
Let $\Pi$ denote the orthogonal projection onto the subspace orthogonal
to the vector $\bf g$, so that
\be\label{def:Pi}
\Pi = I - {\bf g} {\bf g}^{\rm T}
\ee
Then the equation (\ref{start2}) is equivalent to the
following two equations:
\be\label{justZ}
\Pi S  Z \Pi + \Pi Z  S^{\rm T} \Pi = 0
\ee
and
\be\label{Zandw}
S Z {\bf g} + Z {\bf h} + 
\bra {\bf w}, {\bf k-h} \ket {\bf w} = 0
\ee
Furthermore, an application of the Cauchy-Schwarz inequality shows that
$Y \geq {\bf w} {\bf w}^{\rm T}$, and hence
\be\label{Zpos}
Z \geq {\bf w} {\bf w}^{\rm T}
\ee
This means that the pair of
matrices $X' = Z - {\bf w} {\bf w}^{\rm T}$ and $Y' = {\bf w} {\bf w}^{\rm T}$
also satisfy (\ref{thm2:1}). Therefore the existence of any pair $(X,Y)$
which satisfy  (\ref{thm2:1}) implies the existence of a pair
$Z - {\bf w} {\bf w}^{\rm T}$ and ${\bf w} {\bf w}^{\rm T}$ satisfying
(\ref{thm2:1}), where $Z$ and $\bf w$ are related by (\ref{justZ}) and
(\ref{Zandw}). Conversely, if $Z$ and $\bf w$ satisfy (\ref{justZ}),
(\ref{Zandw}) and (\ref{Zpos}), then they provide 
a solution of (\ref{thm2:1}). Therefore we have the following
result which describes the solutions of (\ref{thm2:1}).

\begin{lemma}\label{lemma3}
Suppose that $Z \geq 0$ and $\bf w$ satisfy (\ref{justZ}), (\ref{Zandw})
and (\ref{Zpos}). Then the pair
$X' = Z - {\bf w} {\bf w}^{\rm T}$ and $Y' = {\bf w} {\bf w}^{\rm T}$
satisfy (\ref{thm2:1}). Conversely, suppose $(X,Y)$ satisfy  
(\ref{thm2:1}). Let $Z= X+Y$ and
define $\bf w$ by (\ref{def:w}). Then $Z$ and $\bf w$ satisfy
(\ref{justZ}), (\ref{Zandw}), and (\ref{Zpos}).
\end{lemma}

\medskip
We now return to the equation (\ref{Zandw}) and solve for $\bf w$.
Define
\be\label{def:xi}
{\bf \xi} = S Z {\bf g} + Z {\bf h}
\ee
If $Z$ and $\bf w$ satisfy the equations (\ref{justZ}) and
(\ref{Zandw}), then it must be true that
\be\label{cond:Z1}
\bra {\bf h-k}, {\bf \xi} \ket > 0
\ee
This is a condition on the matrix $Z$. If it is satisfied, then
(\ref{Zandw}) can be solved for ${\bf w}$:
\be\label{solvew}
{\bf w} = \bra {\bf h-k}, {\bf \xi} \ket^{-1/2} \, {\bf \xi}
\ee
The condition $Z \geq {\bf w} {\bf w}^{\rm T}$ is equivalent
to $1 \geq {\bf w}^{\rm T} Z^{-1} {\bf w}$. Defining
\be\label{def:F}
F(Z) = \frac{\bra {\bf \xi}, Z^{-1} {\bf \xi} \ket}
{\bra {\bf h-k}, {\bf \xi} \ket}
\ee
we can combine the two conditions (\ref{cond:Z1}) and (\ref{Zpos}) as
\be\label{cond:Z2}
0 < F(Z) \leq 1
\ee
We can now restate Lemma \ref{lemma3} as follows.

\medskip

\begin{lemma}\label{lemma4}
Suppose that $Z \geq 0$ satisfies (\ref{justZ}) and
(\ref{cond:Z2}). Define $\bf w$ by  (\ref{solvew}). Then the pair
$X' = Z - {\bf w} {\bf w}^{\rm T}$ and $Y' = {\bf w} {\bf w}^{\rm T}$
satisfy (\ref{thm2:1}). Conversely, suppose $(X,Y)$ satisfy  
(\ref{thm2:1}), and let $Z= X+Y$. Then $Z$ satisfies
(\ref{justZ}) and (\ref{cond:Z2}).
\end{lemma}

\medskip
Lemma \ref{lemma4} shows a many-to-one correspondence
between the solutions of (\ref{thm2:1}) and the matrices $Z$ satisfying
(\ref{justZ}) and (\ref{cond:Z2}).
Therefore we have reduced the proof of Theorem \ref{thm2} 
to the following problem: suppose that
there is some matrix satisfying (\ref{justZ}) and
(\ref{cond:Z2}). Then we want to show that there is another such $Z$ 
satisfying (\ref{justZ}) and
(\ref{cond:Z2}) for which both
$X' = Z - {\bf w} {\bf w}^{\rm T}$ and $Y' = {\bf w} {\bf w}^{\rm T}$
are rank 1. Equivalently, we want to show that there is a rank 2 matrix $Z$
satisfying (\ref{justZ}) for which $F(Z)=1$.

\subsection{Representation using Hankel matrices, and the moment problem}
We use the easily verified fact that every symmetric matrix $Z$ which
satisfies (\ref{justZ}) has the following form:
\be\label{Hank}
Z = \pmatrix{z_0 & 0 & -z_1 & 0 & z_2 & \dots \cr
0 & z_1 & 0 & -z_2 & 0 & \dots \cr
-z_1 & 0 & z_2 & 0 & -z_3 & \dots \cr
0 & -z_2 & 0 & z_3 & 0 & \dots \cr
\vdots && \ddots && \dots & \dots \cr}
\ee
Except for the minus signs, (\ref{Hank}) is an example of a Hankel matrix.
Positivity of $Z$ requires that $z_i \geq 0$ for all $i=0,\dots,n-1$, and also imposes other
constraints. To describe the possible values of $\{z_i\}$, 
we introduce for each real $x$ the following  rank 2 matrix of the form (\ref{Hank}):
\be\label{def:Z(x)}
Z(x) = \pmatrix{1 & 0 & -x & 0 & x^2 & \dots \cr
0 & x & 0 & -x^2 & 0 & \dots \cr
-x & 0 & x^2 & 0 & -x^3 & \dots \cr
0 & -x^2 & 0 & x^3 & 0 & \dots \cr
\vdots && \ddots && \dots & \dots \cr}
\ee

\medskip
The next result shows that every matrix $Z \geq 0$ satisfying (\ref{justZ}) and (\ref{Zandw})
can be written as a positive linear combination of the matrices (\ref{def:Z(x)}) for
some non-negative values of $x$.

\medskip
\begin{thm}\label{thm.mom}
Suppose that the $n \times n$ matrix $Z \geq 0$ satisfies equations (\ref{justZ}), (\ref{Zandw})
and (\ref{Zpos}).
Then there is an integer $p \leq (n+1)/2$, non-negative numbers $0 \leq x_0 < x_1 < \dots < x_{p-1}$, and
positive numbers $\mu_0, \dots, \mu_{p-1}$, such that
\be\label{thm.mom.1}
Z = \sum_{i=0}^{p-1} \mu_i Z(x_i)
\ee
If $p=(n+1)/2$ then $x_0 = 0$.
\end{thm}

\medskip
The proof of Theorem \ref{thm.mom} is presented in Appendix A, using standard results
form the theory of moments. Indeed, the result is equivalent to the solution of the
discrete moment problem for $z_0, \dots, z_{n-1}$, which is the problem of finding
distinct non-negative numbers $x_0, \dots, x_{p-1}$ and
positive $\mu_0, \dots, \mu_{p-1}$, so that
\be\label{mom}
z_j = \sum_{i=0}^{p-1} \mu_i x_{i}^j, \quad\quad
0 \leq j \leq n-1
\ee
If there is a solution of (\ref{mom}), then the special form of the Hankel
matrices implies immediately that (\ref{thm.mom.1}) holds, and vice versa.

\medskip
The decomposition (\ref{thm.mom.1}) is the key for solving the problem posed after
Lemma \ref{lemma4}: we will show that if $Z$ satisfies the 
equations (\ref{justZ}), (\ref{Zandw}) and (\ref{Zpos})
so that the representation (\ref{thm.mom.1}) holds,
then at least one of the matrices $\{Z(x_i)\}$ must also satisfy these
equations, and hence by Lemma \ref{lemma4} it provides a solution
of (\ref{thm2:1}). Since $Z(x_i)$ has rank 2, this almost completes the proof
of Theorem \ref{thm2}. The only remaining obstacle is that
the matrix $X' = Z(x_i) - {\bf w} {\bf w}^{\rm T}$ may have rank 2, or equivalently
$F(Z(x_i)) < 1$. To complete the proof we will show that in this case there must be
another number $y < x_i$ for which $F(Z(y)) = 1$. Then
$Z(y) - {\bf w} {\bf w}^{\rm T}$ will have rank 1, and hence can be written 
as ${\bf v} {\bf v}^{\rm T}$, so that (\ref{thm2:2}) holds.

\subsection{Completion of the proof}
The key to the solution is the following Lemma, which displays a remarkable
property of the function $F(Z)$ when $Z$ has the form (\ref{thm.mom.1}).

\begin{lemma}\label{lemma5}
Suppose that $Z$ satisfies the hypotheses of Theorem \ref{thm.mom}, so that
the representation (\ref{thm.mom.1}) holds.
For each $i=0, \dots, p-1$ define
\be\label{def:xi.i}
\xi_{i} = S Z(x_i) {\bf g} + Z(x_i) {\bf h}
\ee
Then
\be\label{lemma5:1}
F(Z) = 
F\left( \sum_{i=0}^{p-1} \mu_i Z(x_i) \right) =
\frac{\sum_{i=0}^{p-1} \, \mu_i \, \bra \xi_{i} , Z(x_i)^{-1} 
\, \xi_{i} \ket}
{\sum_{i=0}^{p-1} \, \mu_i \, \bra {\bf h-k} , \xi_{i} \ket}
\ee
\end{lemma}

\medskip
Lemma \ref{lemma5} in proved in Appendix B. We now use it to complete the
proof of Theorem \ref{thm2}. By assumption there is a matrix $Z$ 
satisfying the hypotheses of Theorem \ref{thm.mom}, and so
Lemma \ref{lemma5} can be used to evaluate $F(Z)$. We now consider how
the right side of (\ref{lemma5:1}) varies as the parameters $\mu_i$ change.
Our goal is to show that there is some $i$ such that
$F(Z) \geq F\left( Z(x_i) \right)$.

Notice first that if $\bra {\bf h-k} , \xi_{i} \ket \leq 0$ for any
$i$, then we do not increase $F(Z)$ by setting $\mu_i=0$, so we will assume 
that $\bra {\bf h-k} , \xi_{i} \ket > 0$ for all $i$.
It is straightforward to compute the derivative with respect to the parameter
$\mu_i$. Since the sign of this derivative is independent of the value of $\mu_i$,
it follows that $F(Z)$ is a monotone function of each $\mu_i$.
Consider first how $F(Z)$ behaves as $\mu_0$ varies.
$F(Z)$ must decrease either as $\mu_0 \rightarrow \infty$
or as $\mu_0 \rightarrow 0$. In the former case we get
\be
F(Z) \geq F\left( Z(x_0) \right)
\ee 
while in the latter case
\be
F(Z) \geq F\left( \sum_{i=1}^{p-1} \mu_i Z(x_i) \right)
\ee
By repeating the same argument if necessary with $\mu_1, \mu_2, \dots$ we eventually
deduce that
\be\label{ineq1}
1 \geq F(Z) \geq F\left( Z(x_i) \right) > 0
\ee 
for some $i$. Recall from Lemma \ref{lemma4} that any matrix
$Z$ satisfying (\ref{justZ}) and $F\left( Z \right) \leq 1$ provides a 
solution of (\ref{thm2:1}), namely $X' = Z - {\bf w} {\bf w}^{\rm T}$
and $Y' = {\bf w} {\bf w}^{\rm T}$. Hence (\ref{ineq1}) implies that
$Z(x_i)$ provides such a solution.

Since $Z(x_i)$ has rank 2, it follows that
$X'$ has rank 1 or 2. If $X'$ has rank 1, then we can
immediately deduce that (\ref{thm2:2}) holds, and the proof is complete. The condition that
$X'$ has rank 1 is $F\left( Z(x_i) \right) = 1$, 
so we are left with the case where $F\left( Z(x_i) \right) < 1$. We now show that
in this case there is another value $y < x_i$ such that
\be\label{def:y}
F(Z(y)) = 1
\ee
This fact follows from these observations:
\begin{itemize}
\item[(a)] $F(Z(x))$ is a rational function of $x$;
\item[(b)] $F(Z(x)) \neq 0$ for all $x \geq 0$;
\item[(c)] either $F(Z(0)) > 1$ or $F(Z(0)) < 0$.
\end{itemize}
To see that (b) is true, note that $F(Z(x)) = 0$ would imply
$S Z(x) {\bf g} + Z(x) {\bf h} = 0$. However this would imply that
$A Z + Z A^{\rm T} = 0$, which is impossible because $A$ is stable.
Similarly (c) is the statement that $Z(0)$ cannot arise as a solution of
(\ref{thm2:1}). Continuity now implies that there must be some $y < x_i$
such that (\ref{def:y}) holds, and this completes the derivation of
(\ref{thm2:2}).

\vskip1in
\noindent{\bf Acknowledgements}
The authors thank R. Shorten
for suggesting the topic of this paper, and for helpful discussions on the
background to the CQLF problem.
This work was supported in part by
National Science Foundation Grant DMS--0101205. 
{~~}

\appendix
\section{Proof of Theorem \ref{thm.mom}}
There are many results known for the discrete moment problem.
We have found special cases of our result in the literature (for example
when $n$ is even and $Z$ is nonsingular \cite{Gant}) but not the full statement.
For this reason we include the proof here;
our starting point is the following result
from the text of Ahiezer and Krein \cite{AK}.

\begin{lemma}\label{lemmaAK}
Given a sequence $s_0, s_1, \dots, s_{2m-2}$ let $K$ denote the
$m \times m$ Hankel matrix with entries
\be\label{def:K}
K_{ij} = s_{i+j-2}, \quad\quad i,j=1,\dots,m
\ee
If $K \geq 0$, then there is an integer $p \leq m-1$, numbers
$x_1, \dots, x_p$, and positive numbers $\mu_1, \dots, \mu_p$
and $M$ such that
\be\label{AKrep}
s_{k} & = & \sum_{i=1}^p \mu_i x_{i}^k \quad\quad
(k=0,\dots,2m-3) \\
s_{2m-2} & = & \sum_{i=1}^p \mu_i x_{i}^{2m-2} + M
\ee
\end{lemma}

\medskip
We now use Lemma \ref{lemmaAK} to prove Theorem \ref{thm.mom}.
First suppose that $n=2m-1$ is odd, and that $Z$ has the form (\ref{Hank}).
Then letting $s_k = z_k$ for $k=0,\dots,2m-2$ it follows that the
matrix $K$ defined by (\ref{def:K}) is positive semidefinite
(the minus signs of some off-diagonal entries in $Z$ can be removed by conjugation
with a diagonal matrix with $\pm 1$ on the diagonal).
Hence the representation (\ref{AKrep}) implies that there is some
$p \leq m-1$ such that
\be\label{outrep}
z_{k} & = & \sum_{i=1}^p \mu_i x_{i}^k \quad\quad
(k=0,\dots,2m-3) \\
z_{2m-2} & = & \sum_{i=1}^p \mu_i x_{i}^{2m-2} + M
\ee
and therefore we get the following representation for $Z$:
\be\label{Zrep1}
Z = \sum_{i=1}^p \mu_i Z(x_{i}) + M {\bf g} {\bf g}^{\rm T}
\ee
The rank 2 matrix $Z(x)$ was defined in (\ref{def:Z(x)}).
It can be written as
\be\label{Zrep2}
Z(x) = {\bf u}(x) {\bf u}(x)^T + x {\bf v}(x) {\bf v}(x)^T
\ee
where ${\bf u}(x)$ and ${\bf v}(x)$ are the $n$-vectors
\be
{\bf u}(x) = \pmatrix{1 \cr 0 \cr -x \cr 0 \cr (-x)^2 \cr
0 \cr (-x)^3 \cr \vdots}, \quad & {\bf v}(x) = S^{\rm T} {\bf u}(x) =
\pmatrix{0 \cr \ 1 \cr 0 \cr -x \cr 0 \cr (-x)^2
\cr 0 \cr \vdots}
\ee
Since $x_{1},\dots,x_p$ are distinct, and $2p \leq 2m-2 < n$,
it follows that the $2p+1$ vectors
$\{{\bf u}(x_i),{\bf v}(x_i)\}$ and $\bf g$ are linearly independent (this is easily
demonstrated using the Hadamard determinant). Hence positivity of $Z$
implies that each term in (\ref{Zrep1}) is separately positive.
From (\ref{Zrep2}) it then follows that 
\be
x_{i} \geq 0, \quad\quad i=1,\dots,p
\ee
It remains to show that $M=0$ in (\ref{Zrep1}). This requires using
(\ref{Zandw}) and (\ref{Zpos}). Suppose first that
$Z$ is singular; then counting dimensions in (\ref{Zrep1}) shows that
either $M=0$ or else $p \leq m-2$. In the latter case it follows that the
vectors $\{{\bf u}(x_i),{\bf v}(x_i)\}$, $\bf g$, $S {\bf g}$ are linearly
independent. Furthermore for any $x$,
\be
S Z(x) {\bf g} = (-x)^{(n+1)/2} {\bf v}(x)
\ee
Therefore from (\ref{Zandw}) it follows that $\bf w$ is a linear combination
of the vectors $\{{\bf u}(x_i),{\bf v}(x_i)\}$, $\bf g$, $S {\bf g}$.
If $M>0$ then the coefficient of $S {\bf g}$ in this combination is nonzero.
But (\ref{Zrep1}) then implies that $\bf w$ is not in the range of $Z$,
which means that (\ref{Zpos}) cannot hold. Therefore we must have $M=0$.

In the general case where $Z$ is non-singular, we argue as follows.
Notice that $Z(0)$ is the matrix whose $(1,1)$ entry is $1$, and all other entries
are zero. If $Z > 0$ then there is $c > 0$ such that
\be
Z = c Z(0) + W
\ee 
where $W \geq 0$ is singular. Since $W$ satisfies (\ref{justZ}),
Lemma \ref{lemmaAK} leads to the representation (\ref{Zrep1}) for $W$,
that is
\be\label{Zrep3}
Z = c Z(0) + \sum_{i=1}^p \mu_i Z(x_{i}) + M {\bf g} {\bf g}^{\rm T}
\ee
where all $x_i >0$.
Since $W$ is singular, either $M=0$ and $p \leq m-1$, or else
$p \leq m-2$. In the former case $2p \leq 2m-2 < n$ so this
establishes (\ref{thm.mom.1}). In the latter case the vectors 
$\{{\bf u}(x_i),{\bf v}(x_i)\}$, $u(0)$, $\bf g$ and $S {\bf g}$ are linearly
independent, and $S {\bf g}$ is not in the range of $Z$. Hence by the same
argument we must have $M=0$. Again we have $2p < n$.

This completes the argument for the case when $n$ is odd.
When $n$ is even, we first create a $(n+1)\times(n+1)$ matrix
${\tilde Z}$ by adding an extra row and column to the matrix $Z$. The new entries
are chosen so that ${\tilde Z}$ has the form (\ref{Hank}). 
This determines uniquely all the entries of ${\tilde Z}$ except
the bottom right corner. This entry is chosen large enough so that ${\tilde Z} \geq 0$.
To see that this is possible, note that the first $n$ entries of the last column of
${\tilde Z}$ are the vector $S Z {\bf g}$. From (\ref{def:F}) we have the bound
\be
|| Z^{-1/2} S Z {\bf g} || & \leq & || Z^{-1/2} (S Z {\bf g} + Z {\bf h}) || 
+ || Z^{1/2} {\bf h} || \nonumber \\
& \leq & F(Z) \bra {\bf h-k}, (S Z {\bf g} + Z {\bf h}) \ket^{1/2} 
+ \bra {\bf h}, Z {\bf h} \ket^{1/2}
\ee
Let $k$ be the bottom right corner entry of ${\tilde Z}$. Then taking
\be
k  \geq \bigg( \bra {\bf h-k}, (S Z {\bf g} + Z {\bf h}) \ket^{1/2} 
+ \bra {\bf h}, Z {\bf h} \ket^{1/2} \bigg)^2
\ee
it follows that ${\tilde Z} \geq 0$. Hence Lemma \ref{lemmaAK} can be applied and we
deduce that
\be\label{Zrep4}
{\tilde Z} = \sum_{i=1}^p \mu_i {\tilde Z}(x_{i}) + M {\bf {\tilde g}} {\bf {\tilde g}}^{\rm T}
\ee
where ${\tilde Z}(x_{i})$ and ${\bf {\tilde g}}$ are the $(n+1)$ 
dimensional versions, and again $x_{i} \geq 0$, and 
$\mu_i > 0$. Since $2p \leq n$ we immediately deduce
(\ref{thm.mom.1}) by restricting both sides of (\ref{Zrep4}) to the top left $n \times n$ block.
This completes the proof of Theorem \ref{thm.mom}.

\section{Proof of Lemma \ref{lemma5}}
We assume that the representation (\ref{thm.mom.1}) holds. 
Lemma \ref{lemma5} follows from a result in linear algebra which
we state and prove in Lemma \ref{inverselemma} below. 
First we verify that the conditions of the lemma are satisfied.

From (\ref{Zrep2}) it follows that the range of $Z$ is spanned by the
vectors $\{{\bf u}(x_i), {\bf v}(x_i)\}$ ($i=0,\dots,p-1$).
If $2p \leq n$ then these vectors are independent, and hence 
$rk(Z)$ = $\sum rk(Z(x_i))$. If $2p = n+1$, then $n$ is odd, and 
Theorem \ref{thm.mom} implies that $x_0 = 0$. Since the vectors
${\bf u}(0)$ and $\{{\bf u}(x_i), {\bf v}(x_i)\}$ ($i=1,\dots,p-1$) 
are independent, it is again true that 
$rk(Z)$ = $\sum rk(Z(x_i))$. 

Furthermore, 
\be
S Z(x_i) {\bf g} &=& \left \lbrace \begin{array}{ll} -
(-x_i)^{n/2} {\bf u}(x_i) & n \mbox{ is even} \\
\\
(-x_i)^{(n+1)/2} {\bf v}(x_i)
&  n \mbox{ is odd} \end{array} \right. \\
&& \nonumber \\
& \in & Ran(Z(x_i))
\ee
Since $Z(x_i) {\bf h}$ is clearly in $Ran(Z(x_i))$, we deduce that
\be
 \xi_i = Z(x_i) {\bf h} + S Z(x_i) {\bf g}
 \,\, \in \, Ran(Z(x_i))
\ee

\medskip
\begin{lemma}\label{inverselemma}
Let $Z$ be an $n \times n$ matrix
\be
Z = \sum_i \mu_i Z_i
 \ee
such that each $Z_i$ is symmetric and rk($Z$) = $\sum
rk(Z_i)$.
Also, let 
\be
{\bf v} = \sum_i \mu_i {\bf v}_i
\ee
where each ${\bf v}_i \in Ran(Z_i)$. Then:
\be
\bra {\bf v}, Z^{-1} {\bf v} \ket = \sum_i \mu_i \bra {\bf v}_i,
Z_i^{-1} {\bf v}_i \ket
\ee
\end{lemma}

\medskip

\pf
It is not assumed that $Z$ is invertible on $R^n$;
since ${\bf v} \in Ran(Z)$, $Z^{-1} {\bf v} $ is always
well-defined. We first prove the result in the case that each $Z_i =
\lambda_i {\bf u}_i {\bf u}_i^T$ is rank 1. Since each ${\bf v}_i$ is in
the range of $Z_i$, 
\be
{\bf v}_i  =  a_i {\bf u}_i, \quad\quad
{\bf v}  =  \sum_i \mu_i a_i {\bf u}_i \label{defv}
\ee
The fact that  rank($Z$) = $\sum \mbox{ rank}(Z_i)$
implies that the ${\bf u}_i$'s form a basis for Ran(Z).
We write $Z^{-1} v$ in that basis with arbitary
coefficients:
\be
Z^{-1} {\bf v} = \sum_j \alpha_j {\bf u}_j
\ee
Now apply $Z$ to both sides:
\be
{\bf v}  = Z(\sum_j \alpha_j {\bf u}_j) 
 =  \sum_{i,j} \mu_i Z_i \alpha_j {\bf u}_j
 = \sum_{i,j} \mu_i \lambda_i  \bra {\bf u}_i, {\bf u}_j \ket
\alpha_j {\bf u}_i \label{newdefv}
\ee
Comparing coefficients in (\ref{defv}) and
(\ref{newdefv}), we see that for all $i$,
\be
a_i =  \lambda_i  \sum_j \bra {\bf u}_i, {\bf u}_j \ket \alpha_j
\label{key}
\ee
Now we can calculate:
\be
\bra {\bf v}, Z^{-1} {\bf v} \ket & = & \sum_{i,j} \mu_i a_i \bra
{\bf u}_i, {\bf u}_j \ket \alpha_j \\
& = & \sum_{i} \mu_i a_i
\left(\frac{a_i}{\lambda_i}\right) \\
& = & \sum_{i} \mu_i 
\left(\frac{a_i^2}{\lambda_i}\right)\\
& = & \sum_{i} \mu_i \bra {\bf v}_i, Z_i^{-1} {\bf v}_i \ket
\ee

To get the full result is now straightforward: since each
$Z_i$ is symmetric, it can written in terms of an
orthonormal basis:
\be
Z_i = \sum_{k = 1}^{r_i} \lambda_{i,k} {\bf u}_{i,k} {\bf u}_{i,k}^{\rm T}, \quad\quad
v_i = \sum_{k = 1}^{r_i} a_{i,k} {\bf u}_{i,k}
\ee
Now, we can write $Z$ as a sum of linearly independent
rank 1 projections and apply what was shown above:
\be
 Z  & = &  \sum_i \mu_i Z_i
  =  \sum_{i,k} \mu_i \lambda_{i,k} {\bf u}_{i,k}
{\bf u}_{i,k}^T \\
 {\bf v} & = & \sum_i \mu_i {\bf v}_i
  =  \sum_{i,k} \mu_i a_{i,k} {\bf u}_{i,k} \\
 \bra {\bf v}, Z^{-1} {\bf v} \ket & = & \sum_{i,k} \mu_i 
\left(\frac{a_{i,k}^2}{\lambda_{i,k}}\right)
\label{end1}
 \ee
The final observation is that, for each $i$, the
${\bf u}_{i,k}$ are orthogonal, which means  
  \be
  Z_i^{-1}& =& \sum_{k = 1}^{r_i}
\left(\frac{1}{\lambda_{i,k}}\right) {\bf u}_{i,k} {\bf u}_{i,k}^T
\\
 \bra {\bf v}_i, Z_i^{-1} {\bf v}_i \ket &=& \sum_{k = 1}^{r_i}
\left(\frac{a_{i,k}^2}{\lambda_{i,k}}\right)
\label{end2}
\ee
Combining (\ref{end1}) and (\ref{end2}), we see
 \be
  \bra {\bf v}, Z^{-1} {\bf v} \ket  =  \sum_{i,k} \mu_i 
\left(\frac{a_{i,k}^2}{\lambda_{i,k}}\right) 
   =  \sum_i \mu_i \bra {\bf v}_i, Z_i^{-1} {\bf v}_i \ket
  \ee
which was to be shown.

\end{document}